\documentclass[]{article}
\usepackage[left=2cm,right=2cm,top=2.54cm,bottom=2.54cm]{geometry}
\usepackage[utf8]{inputenc}
\usepackage{graphicx}
\usepackage{verbatim}
\usepackage[numbers]{natbib}
\usepackage{adjustbox}
\usepackage{bm}
\usepackage{hyperref}
\usepackage[many]{tcolorbox}
\usepackage{amssymb,amsthm,amsmath,amsfonts}
\usepackage{tikz}
\usepackage{authblk}
\usetikzlibrary{shapes,arrows,positioning}

\usepackage[ruled,vlined]{algorithm2e}
\SetKwInput{KwInput}{Input}                
\SetKwInput{KwOutput}{Output}              

\theoremstyle{plain}
\newtheorem{thm}{Theorem}[section]

\theoremstyle{definition}

\theoremstyle{remark}
\newtheorem*{rem}{Remark}

\tikzset{block/.style={rectangle, draw, fill=gray!15,
		text width=16em, text centered, rounded corners, minimum height=4em},
	line/.style={draw, -latex'}}

	\title{Automated solving of constant-coefficients second-order
	linear PDEs using Fourier analysis}
\author[1]{E. Roque}
\author[2]{J. A. Vallejo}
\affil[1]{Departamento de Matemáticas\\ 
	Cinvestav Unidad Querétaro\\ 
	Libramiento Norponiente \#2000, Fracc. Real de Juriquilla\\
	Quer\'{e}taro, Qro., 76230 México,
	\texttt{earoque@math.cinvestav.mx}}
\affil[2]{Departamento de Matemáticas Fundamentales\\ 
	Universidad Nacional de Educación a Distancia\\ 
	C. Juan del Rosal 10, CP 28040, Madrid, Spain\\
	\texttt{jvallejo@mat.uned.es}}
\begin{document}
	\maketitle
		\begin{abstract}
		We provide the details of an implementation of Fourier techniques for solving second-order
		linear partial differential equations (with constant coefficients) using a computer algebra
		system. The general Sturm-Liouville problem for the heat, wave and Laplace operators on
		the most common bounded domains is covered, as well as the general second-order linear
		parabolic equation with constant coefficients, which includes cases such as the
		convection-diffusion equation, by reduction to the heat equation.
	\end{abstract}
	

	\section{Introduction}
		The need for explicit solutions of differential equations (be they ordinary or partial) arises in the everyday practice of physicists, engineers, and scientists in general. The class of ordinary differential equations that can be solved by elementary methods, mainly clever changes of variables and linear algebra, is quite ample. Accordingly, there exists a vast amount of software devoted to finding both their explicit and numerical solutions. This is in contrast with the case of partial differential equations (PDEs). Here, solution techniques are more advanced, involving Fourier analysis for bounded domains and integral transforms in the unbounded case\footnote{In this paper we restrict our attention to classical solutions, thus leaving aside considerations regarding weak solutions and the use of such tools as distribution theory.}. Although there are a good deal of software packages related to their numerical computations (most  of them using the finite element method), there is a shortage of options for finding explicit analytic solutions. In fact, to the best of our knowledge,
		only Maplesoft's Maple\texttrademark\  and Wolfram's Mathematica\texttrademark\ --both closed and proprietary software-- are
		capable of this task in a more or less general setting\footnote{SymPy (a Python library for symbolic mathematics) offers basic support for first-order PDEs.}. Both are excellent general-purpose computer algebra systems (CAS) and do their job fast and efficiently, so one could wonder what the motivation is for considering the problem of automated solution of PDEs or thinking it is not already solved. To this question, we can offer the following answers:
	\begin{enumerate}
		\item Although they are excellent at finding solutions, still there are examples where Maple\texttrademark\
		or Mathematica\texttrademark\  cannot find it (see examples in subsection \ref{secex}).
		In those cases, it is very difficult
		to determine the origin of the failure, as both are closed-source packages.
		The user can have a hard time trying
		to discern whether there is an issue in the formulation of the problem or in the method of solution.
		\item Proprietary software usually has a cost, and this case is no exception. While both packages are worth
		their price, software license costs can be a serious hurdle for people working in underdeveloped countries, students in general,
		or even freelancers for whom an open-source alternative would be desirable.
		\item From a pedagogical point of view, the use of closed-source software is like using a black box for
		finding answers. The teacher can expose an algorithm in the classroom and then say it is exemplified by the
		output, but really there is no clue about what is being actually computed: if there is any mathematical trick used that is not contemplated in the algorithm, or if use is made of pre-computed tables of particular cases, and so on.
		\item Lastly, the problem is interesting in itself from a mathematical and computational perspective
		due to the blending of purely symbolic and numeric algorithms it requires, as will be seen in detail below.
	\end{enumerate}
	
	These reasons led us to write a library, called \texttt{pdefourier}, for solving the general
	Sturm-Liouville problem for the heat, wave and Laplace operators with sources, as well as some particular
	problems reducible to these, such as the convection-diffusion equation.
	Our implementation uses the open-source Maxima CAS \cite{Max},
	which has a programming language very closely tied to LISP and should be translatable to any other CAS or
	programming language in general, such as SymPy or Symbolic C++, without difficulty. \\
	
	Let us summarize the contents of the paper. Section \ref{sec-slp} recalls the main
	theoretical results on which the paper is based. In section \ref{sec-hurdles} we comment on the challenges posed by the automated computation of solutions using
	Fourier techniques in general. By including examples of computations using Mathematica\texttrademark,
	we present a panorama of the current state of the art.\footnote{We must insist in this aspect, since nothing could
		be further from our minds than to do a comparison test or a benchmark against a commercial software, supported by
		a whole team of scientists.}
	Section \ref{sec-fcoeff} describes how we handle the issues mentioned in section \ref{sec-hurdles}, detailing some of the
	core features of the library, and we finish in section \ref{sec-examples} by showing some examples and applications.
	
	\begin{rem}
		In what follows, all Mathematica\texttrademark\ commands have been executed using version
		14.2.1.0 on a Windows platform. The Maxima commands have been executed on the same Windows machine,
		with version 5.47.0 compiled against the LISP implementation SBCL 2.3.2.
		It is assumed that the package \texttt{pdefourier} has already been loaded. The first time the user loads the package \texttt{pdefourier}, compilation messages might get printed on the screen since the package loads \texttt{lapack}, which needs to be compiled the first time. 
		For installation instructions see Appendix \ref{app}. Maple\texttrademark\ commands were executed using the version 2026 in a Windows machine.
	\end{rem}

	\section{Sturm-Liouville problems and eigenfunctions expansion}\label{sec-slp}
	
	From a mathematical point of view, the basic models given by the heat, wave and
	Laplace equations are all second-order linear partial differential equations with constant
	coefficients. They are supplemented by some initial and/or boundary conditions (the general situation contemplates both, and it is known as an IBVP or initial-boundary value problem), and in this
	work we will assume that the problem is set on a bounded spatial domain $[0,L]\subset\mathbb{R}$,
	so these conditions will be given at the boundary points $x=0,x=L$ (the only exception
	being Laplace's equation, which will be considered on a plane domain).
	The general problem for the heat and wave equations can be written as:
	\begin{equation}\label{problem}
		\begin{cases}
			\mathcal{D}u =Q(x,t) \\
			u(x,0)=F(x) \\
			u_x(x,0)=G(x) \\
			\alpha_1 u(0,t) +\beta_1 u_x(0,t)= h_1(t) \\
			\alpha_2 u(L,t) +\beta_2 u_x(L,t)= h_2(t)\,,
		\end{cases}
	\end{equation}
	where $\mathcal{D}$ denotes either the heat or the wave operator. In the case of the Laplace's
	equation we consider a wider class of plane domains (rectangles, disks, wedges, annuli) and boundary conditions. In a rectangular domain \([0,a] \times [0,b] \) we consider the problem with mixed boundary conditions of the form:
	\begin{equation}\label{eq-lprect}
		\begin{cases}
			\Delta u=0 \quad \alpha, \beta, \gamma, \delta \in \{ 0,1 \} \\
			(1-\alpha)u(x,0)+\alpha u_y(x,0)=f_0(x) \\
			(1-\beta)u(x,b)+\beta u_y(x,b)=f_b(x) \\
			(1-\gamma) u(0,y) + \gamma u_x (0,y)=g_0(y) \\
			(1-\delta) u(a,y) + \delta u_x (a,y)= g_a(y)
		\end{cases}
	\end{equation}
	
	On the remaining domains, the Laplace equation is solved in polar coordinates with Dirichlet and Neumann boundary conditions; the only exception being the case of an annular domain, where only Dirichlet boundary conditions are considered. Here is a
	summary of the cases:
	
	\begin{equation}\label{eq-lpdwa}
		\begin{cases}
			u_{rr}+\frac{1}{r}u_r+\frac{1}{r^2}u_{\theta\theta}=0 \\
			(r,\theta) \in (0,R_1]\times [0,\alpha[, \, 0<\alpha\leq 2\pi, \,0<R_1 \text{ (Disk/Wedge)} \\
			(r,\theta) \in [R_1,R_2]\times [0,2\pi[, \, 0<R_1<R_2 \text{ (Annulus)} \\
			\text{Dirichlet boundary conditions:} \\
			u(r,0)=0=u(r,\alpha) \text{ if } \alpha<2\pi \text{ (Wedge)} \\
			u(R_1,\theta)=f(\theta)  \text{ (Annulus/Disk/Wedge)}\\
			u(R_2,\theta)=g(\theta) \text{ (Annulus)} \\
			\text{Neumann boundary conditions:} \\
			u(r,0)=0=u(r,\alpha) \text{ if } \alpha<2\pi \text{ (Wedge)}\\
			u_r(R_1,\theta)=f(\theta) \text{ (Disk/Wedge)}\\
		\end{cases}
	\end{equation}
	
	A common procedure to solve these equations is to apply separation of variables: Writing
	$$
	u(x,t)=X(x)T(t)\,,
	$$
	and denoting generically $X$ or $T$ by $\phi$, leads to a regular Sturm-Liouville problem
	whose general form is
	\begin{equation}\label{SL}
		\mathcal{L}\phi := (p\cdot\phi')'+(q +\lambda r)\phi =0\,,
	\end{equation}
	where $p,r>0$ are positive real functions, and $p',q,r$ are continuous on $[0,L]$. The abstract theory of Sturm-Liouville operators is very well developed (see \cite{AG08,Hel67}), and its main result can be summarized as follows:
	
	\begin{thm}
		Let $\mathcal{L}$ be a Sturm-Liouville operator as in \eqref{SL}, subject to the boundary conditions
		derived from those in \eqref{problem} in the homogeneous case. Assume that, for $i\in\{1,2\}$,
		\begin{equation}\label{conds}
			|\alpha_i|^2+|\beta_i|^2\neq 0\,.
		\end{equation}
		Then, the following hold:
		\begin{enumerate}
			\item All the eigenvalues $\lambda$ are real and form a denumerable set.
			\item The set of all eigenfunctions, $\mathcal{E}$, is an orthogonal set.
			\item The set $\mathcal{E}$ is complete in the space of square-integrable functions
			$L^2([0,L])$.
			\item The sequence $(\lambda_n)$ diverges to $+\infty$ and is bounded from below.
			\item If $\alpha_i,\beta_i,p,-q$ are all non-negative, then the eigenvalues $\lambda_n$ are
			all non-negative.
		\end{enumerate}
	\end{thm}
	
	Notice that the preceding theorem applies to \emph{homogeneous} initial and boundary conditions and \(Q=0 \).
	To deal with the non-homogeneous case , we use the method of eigenfunction expansion or finite cosine or sine transform (see \cite{Hab13,Wal98}).
	Once we know how to deal with the general Sturm-Liouville problem for the heat, wave and Laplace
	operators, we can solve some other problems that can be reduced to one of these. The cases of the
	convection-difussion equation and heat loss through lateral boundaries are then particular cases.
	Thus, consider the general expression of a second-order linear parabolic equation with constant
	coefficients with sources:
	\begin{equation}\label{auxiliar}
		u_t=\kappa u_{xx} + v u_x + cu +Q(x,t), 
	\end{equation}
	Let us apply the change of variable
	\[
	u=we^{\rho(x,t)}\,,
	\]
	where
	\[
	\rho(x,t)=(c-\frac{v^2}{4\kappa})t-\frac{v}{2\kappa}x
	\]
	Substituting into \eqref{auxiliar} we get
	\begin{align*}
		Q&=u_t-\kappa u_{xx} -vu_x-cu \\
		&=e^\rho(w_t-\kappa w_{xx})
	\end{align*}
	hence, the original equation in terms of the new dependent variable reduces to the heat equation with sources
	\[
	w_t=\kappa w_{xx}+e^{-\rho}Q
	\]
	
	Of course, the initial and boundary conditions for $u(x,t)$ must be transformed accordingly into conditions
	for $w(x,t)$. \\
	
	The Sturm-Liouville theory is the basic tool used in our Maxima package \texttt{pdefourier}. The algorithm
	can be summarized as follows:  we first
	solve the eigenvalue problem associated to \eqref{SL}, normalize the corresponding eigenfunctions and then proceed to develop the solution to \eqref{problem} in term of eigenfunctions, taking non-homogeneous conditions into account with the use of the methods commented above. These steps involve some delicate issues that we discuss next.

\section{Difficulties in the computation of Fourier coefficients}\label{sec-hurdles}
	\subsection{Singular values of the Fourier coefficients}\label{ssec:sv}
	Given a function \(f\) with period \(2L\), the Fourier coefficients of \(f\) are
	\begin{equation}\label{eqn:fourier-coeff}
		a_0=\frac{1}{2L}\int^L_{-L}f(x)\,\mathrm{d}x\,,\
		a_n=\frac{1}{L}\int^L_{-L}f(x)\cos \left(\frac{n \pi x}{L}\right)\,\mathrm{d}x\,,\
		b_n=\frac{1}{L}\int^L_{-L}f(x)\sin \left(\frac{n \pi x}{L}\right)\,\mathrm{d}x,\; n\in \mathbb{N}.
	\end{equation}
	It is conventional to write 
	\begin{equation}\label{eqn:fourier-series}
		f \sim a_0+ \sum_{n=1}^{\infty} \left[a_n\cos \left(\frac{n \pi x}{L}\right)+b_n\sin \left(\frac{n \pi x}{L}\right)\right].
	\end{equation}
	The formal series on the right is called the Fourier series of the function \(f\). \\
	
	The computation of Fourier coefficients in a CAS can be done using integration formulas. However, particular care must be taken to separate ``singular values'' that lead to divergences in the general expressions of the Fourier coefficients. \\
	
	There are two types of special input functions whose Fourier coefficients must be calculated separately, namely, the functions belonging to one of the cosine, sine, trigonometric or complex orthogonal systems or a product of a polynomial with them. In these cases, the Fourier coefficients have a general expression
	\begin{equation}\label{eqn:scoeff}
	\frac{C}{\rho (n)}
	\end{equation}
	with \(\rho\) a polynomial, which is only valid for the values of \(n\) where \(\rho(n) \neq 0\). Throughout this work, the set of integers \(n\) for which Fourier coefficients of the form \eqref{eqn:scoeff} satisfy   \(\rho(n)=0\) are called singular values. The process to detect singular values in a given function lies at the bottom of the difficulties experienced by such mature software as Mathematica\texttrademark\ . Due to the
	widespread of this CAS, we will use it to get a benchmark for our library performance. \\
	
	To get a sense of the kind of problems we are talking about, consider for simplicity the case of a \(2 \pi\)-periodic function. In practice, the explicit computations of these coefficients depends heavily on the orthogonality
	relations ($m,n\geq 1$)
	\[
	\frac{1}{\pi}\int^\pi_{-\pi}\cos (mx)\cos (nx)\,\mathrm{d}x =\delta_{mn}\,,\
	\frac{1}{\pi}\int^\pi_{-\pi}\sin (mx)\sin (nx)\,\mathrm{d}x =\delta_{mn}\,.
	\]
	
	However, most computer algebra systems fail to recognize the possibility that $m=n$. Thus,
	Mathematica\texttrademark\ gives
	\begin{verbatim}
		In[1]:= Integrate[Cos[m x] Cos[n x], {x, -Pi, Pi}]
		Out[1]= (2 m Cos[n \[Pi]] Sin[m \[Pi]] - 2 n Cos[m \[Pi]] Sin[n \[Pi]])/(m^2 - n^2)
	\end{verbatim}
	and the same is true of Maxima:
	\begin{verbatim}
		(%i1)	integrate(cos(m*x)*cos(n*x),x,-%pi,%pi);
		(%o1)	(2*((n-m)*sin(%pi*n+%pi*m)+(n+m)*sin(%pi*n-%pi*m)))/(2*n^2-2*m^2)
	\end{verbatim}
	
	\begin{rem}
		At the time we developed the package, in Maxima 5.47 or previous versions, the output for such integral was as above. In the newest Maxima version, which is Maxima 5.49, now you get the following:
		\begin{verbatim}
			(%i1)	declare(n,integer)$
			(%i2)	integrate(cos(m*x)*cos(n*x),x,-%pi,%pi);
			"Is "n" equal to "m"?"y;
			(%o2)	%pi
			(%i3)	integrate(cos(m*x)*cos(n*x),x,-%pi,%pi);
			"Is "n" equal to "m"?"n;
			(%o3)	(2*((n-m)*sin(%pi*n+%pi*m)+(n+m)*sin(%pi*n-%pi*m)))/(2*n^2-2*m^2)
		\end{verbatim}
		Nonetheless, this is still not enough to overcome all the difficulties we explain throughout this work, and it still fails in the rest of scenarios pointed out here unless we note otherwise. The functionality of our package is not affected by the update.
	\end{rem}
	This integration issue has many consequences due to the presence of singular values of $n$, as pointed above.
	Take as an example the computation of Fourier coefficients for
	\begin{equation}\label{fexample}
		f(x)=3x^2\cos(7x)\,.
	\end{equation}
	Using the trigonometric identity
	\begin{equation}\label{eqn:trigid-cos}
	\cos\theta\cos\phi =\frac{1}{2}\cos(\theta-\phi)+\frac{1}{2}\cos(\theta +\phi)
	\end{equation}
	and integrating by parts, it is straightforward to obtain
	\begin{equation}\label{eqn:x2cos7x}
		\begin{split}
		a_0=& -\frac{6}{49} \\
		a_n=& \begin{cases}
			-\dfrac{12(n^2+49)(-1)^n}{(n-7)^2(n+7)^2}\mbox{ if }n\neq 7 \\[7pt]
			\dfrac{98\pi^2+3}{98}\mbox{ if }n=7
		\end{cases}\\
		b_n=& 0
	\end{split}
	\end{equation}
	
	Mathematica\texttrademark\  fails to compute correctly the Fourier coefficients directly from the integral formulas because it does not recognize the singular value corresponding to $n=7$:
	\begin{verbatim}
		In[2]:= Assuming[Element[n,Integers],1/Pi Integrate[Cos[n x] 3 x^2 Cos[7 x],{x,-Pi,Pi}]]
		Out[2]= -((12 (-1)^n (49 + n^2))/(-49 + n^2)^2)
	\end{verbatim}
	and neither does Maxima (not even in version 5.49):
	\begin{verbatim}
		(%i2)	block([n],declare(n,integer),
		factor(integrate(cos(n*x)*3*x^2*cos(7*x),x,-%pi,%pi)/%pi));
		(%o2)	-(12*(n^2+49)*(-1)^n)/((n-7)^2*(n+7)^2)
	\end{verbatim}
	
	\subsection{Equivalence of trigonometric expressions} \label{ssec:equivtrig}
	Another source of concern when dealing with trigonometric expressions is the existence of multiple equivalent ways of writing
	them. This is a well-known issue in different computer algebra systems.
	The set of heuristics rules Mathematica\texttrademark\, appears to use inside its Fourier coefficients built-in functions try to support common special cases, but are highly sensitive to the way the input is written and do not take advantage of the linearity of the integral nor distribute over subintervals of piecewise-defined functions. For instance, using equation \eqref{eqn:trigid-cos} in the case $\theta=\varphi$ we get:
	\begin{equation}\label{trigexample}
		\cos (\theta)^2 = \frac{1}{2}+\frac{1}{2}\cos{\left( 2 \theta \right) }
	\end{equation}
	
	For illustration purposes, here are the results given by Mathematica\textsuperscript{TM} for the Fourier cosine coefficients of each side of the equation.
	\begin{verbatim}
		In[3]:=FourierCosCoefficient[Cos[x]^2, x, n]
		Out[3]:=0
		In[4]:=FourierCosCoefficient[(1+Cos[2x])/2, x, n]
		Out[4]:=1/2 (DiscreteDelta[-2 + n] + 2 DiscreteDelta[n])
	\end{verbatim}
	
	Moreover, its routines do not take advantage of linearity, as a slight change in the input function prevents them from
	working:
	\begin{verbatim}
		In[3]:=FourierCosCoefficient[(2+Cos[2 x])/2,x,n]
		Out[3]:=0
	\end{verbatim}
	
	The strategy followed in our implementation of the package is to internally transform any trigonometric
	function using Maxima's command \texttt{trigrat}, whose output we refer as canonical form in this paper. Please refer to the Maxima's documentation \cite{Max} for further details about \texttt{trigrat}. In particular, powers of sine and cosine are transformed into a linear combination of sines and cosines, in a similar fashion to eq. \eqref{trigexample}. Thus, it becomes easy to decide whether or not the input contains an expression whose Fourier coefficients have singular values using the pattern matching capabilities of the
	Maxima CAS, as described in section \ref{sec-fcoeff}.
	
	\subsection{Piecewise-defined functions}\label{ssec:hpw}
	Let us point that, in applications, many of the functions we deal with are
	piecewise defined and this poses its own challenges. Consider the complex Fourier coefficients of the function \(\sin 3x.\) In Mathematica we can compute them with the following command.
	\begin{verbatim}
		In[5]:= FourierCoefficient[Sin[3x],x,n]
		Out[5]=
	\end{verbatim}
	\[
\begin{cases}
	-i/2 \; & n=3\\
	i/2 \;& n=-3 \\
	0 \;&\text{ True }
\end{cases}
	\]
	 Even though Mathematica is able to detect that \(n=\pm 3\) are special cases, it cannot longer detect that the coefficients with indices \(\pm 3\) need to be computed separately if the expression \(\sin 3x\) appears inside a piecewise-defined function.
	 \begin{verbatim}
	 	In[6]:= f[x_] := Piecewise[{{0, -Pi <= x < 0}, {Sin[3 x], 0 <= x < Pi}}]
	 	In[7]:= FourierCoefficient[f[x],x,n]
	 	Out[7]=
	 \end{verbatim}
	\[
	-\frac{3 \left((-1)^n+1\right)}{2 \pi  \left(n^2-9\right)}
	\]
	Maxima currently does not have built-in support for piecewise-defined expressions which are of particular importance in engineering and physics applications. We developed another package called \texttt{piecewise} automatically loaded by \texttt{pdefourier} to deal with this kind of functions.  In particular, we are able to detect the parity of a piecewise-defined function defined on an interval \([-L,L] \) by the means of comparison of subintervals allowing us to simplify the computation of Fourier coefficients for odd or even functions. Also, we take advantage of the linearity and interval addition properties of the integral; because of this, our pattern matching rules work even with piecewise-defined functions and equivalent trigonometric expressions as we will see in section \ref{sec-fcoeff}.

	\section{Computation of the Fourier coefficients and series}\label{sec-fcoeff}
	
	The aim of this section is to explain the strategies followed in the implementation of our package to tackle the different challenges we have shown in the previous section.
	
	\subsection{Conversion of piecewise-defined functions}
	
	Inherited from Lisp, Maxima's main data structures are lists, thus the natural way to proceed was to translate a piecewise-defined function to a list and viceversa. Having a way to convert a piecewise expression to a list, it is straightforward to write some other functions to make operations between them or to compute the derivative and the integral of a piecewise expression written as a list. Although the package \texttt{piecewise} is able to work in more general settings, we restrict ourselves here to show only how the package works in a bounded interval. The remaining details are available in the package
	documentation. For illustration purposes, consider an if-else expression of the form: \\
	
	\texttt{if x>=a$_0$ and x<=a$_1$ then expr$_1$ elseif ... elseif x>a$_{n-1}$ and x<=a$_{n}$ then expr$_{n}$ } \\
	
	We impose the restrictions $a_i\leq a_{i+1}, a_i \in \mathbb{R},$ to avoid potential issues with the detection of valid intervals in all the equivalent ways of writing them, and to do it efficiently.
	The output returned by \texttt{pw2list} has the form:
	\[
	[[[a_0,a_1],\text{expr}_1],[[a_1,a_2],\text{expr}_2,]\ldots]
	\]
	
	The approach of internally working with lists instead of if-else expressions is particularly handy since the Maxima
	language supports the functional programming paradigm so once we have proper detection of special cases with pattern-matching,
	we can easily map the routines to each subinterval.
	
	\subsection{Computation of Fourier coefficients}

	As mentioned in the preceding section, the product of a polynomial by a trigonometric function gives rise to
	singular values when determining Fourier coefficients, having as particular cases trigonometric functions whose
	wave number is a multiple of $\pi/L$. Because of the linearity of the integral, using trigonometric canonical forms (see subsection \ref{ssec:equivtrig})
	it is enough to detect the following patterns:
	\begin{equation}\label{eqn:patterns}
	x^r\cos \left(\frac{m \pi x}{L}\right), x^r\sin \left(\frac{m \pi x}{L}\right), \cos \left(\frac{m \pi x}{L}\right),\sin \left(\frac{m \pi x}{L}\right), \quad r,m \in \mathbb{N}
	\end{equation}
		To do so, we used the Maxima built-in commands \texttt{defmatch} and \texttt{matchdeclare}. Then, two similar strategies were followed depending on the input expression being piecewise-defined or not. Algorithm \ref{algorithm:fcoeffs} describes the steps followed in our implementation to calculate the Fourier coefficients.
	
	\begin{algorithm}[htp]\label{algorithm:fcoeffs}
		\SetAlgoLined
		\DontPrintSemicolon
		\KwInput{ expression \emph{expr}; variable \emph{var}, semi-length of interval \emph{L} }
		\KwOutput{ A list of Fourier coefficients of the form \([[a_0,a_n,b_n], l.s.v ]\) \\ l.s.v=list of singular values}
		\eIf{\emph{expr} is not piecewise-defined}{
			Apply the simplification command \texttt{trigrat} to convert trigonometric expressions appearing in \emph{expr} into their canonical form.\;
			Expand \emph{expr} fully and convert the expanded expression \( term_1+\ldots+term_s\) into a list \([\text{term}_1,\ldots,\text{term}_s ] \). \;
			If any parity is detected by heuristic routines, use it to simplify computation. In general, map the built-in integration command to each element of the list to compute the general coefficients \( a_n(\text{term}_k), \, b_n(\text{term}_k) \).   \;
			Apply an auxiliary function that searches for the patterns appearing in eq. \eqref{eqn:patterns} to obtain the set of unique indices having singular values, if any, of \( a_n(\text{term}_k), b_n(\text{term}_k) \) and store it in a list \(A\).\;
			Sum over \( k\) the general expression of \( a_n(\text{term}_k), \, b_n(\text{term}_k) \) to obtain the final answer for \(a_n, b_n \) and compute separately the Fourier coefficients of the indices appearing in \(A\) storing the result in the list of singular values. The list of singular values, if any, is written as follows \( [[j,a_j,b_j],\ldots] \), otherwise, an empty list [\,] is returned. \;
			}{
			Apply the above steps to each subinterval of the domain of the piecewise defined \emph{expr}, and sum over \(k^\prime\) both \( a_n(\text{subinterval}_{k^\prime}), \, b_n(\text{subinterval}_{k^\prime}) \) to obtain \(a_n, b_n \).\; 
			Join the sets of indices having singular values in each subinterval, and compute separately the Fourier coefficients of these indices on the whole interval. The list of singular values, if any, is written as follows \( [[j,a_j,b_j],\ldots] \), otherwise, an empty list [\,] is returned. \;
			}
		\caption{Computation of Fourier coefficients}
	\end{algorithm}

	The same ideas also apply for the case of the complex, sine or cosine coefficients; the only difference is in the way the output is returned to the user (see Table \ref{table:formats}).
	
	\begin{table}
		\centering
		\begin{tabular}{llll}
			Type & Command & Answer format & List of s.v format \\
			\hline\\\\[-1.9\medskipamount]
			Trigonometric & \texttt{fouriercoeff(expr, var, L)} & $[[a_0, a_n, b_n], \mbox{list of s.v}]$ & $[\,] \mbox{ / } [[j,a_j,b_j],\ldots]$ \\
			Complex &  \texttt{cfouriercoeff(expr, var, L)} & $[[c_0,c_n], \mbox{list of s.v}]$ & $[\,] \mbox{ / } [[j,c_j],\ldots]$ \\
			Cosine & \texttt{fouriercoscoeff(expr, var, L)} & $[[a_0,a_n], \mbox{list of s.v}]$ & $[\,] \mbox{ / } [[j,a_j],\ldots]$ \\
			Sine & \texttt{fouriersincoeff(expr, var, L)} & $[[b_n], \mbox{list of s.v}]$ & $[\,] \mbox{ / } [[j,b_j],\ldots]$ \\[7pt]
			\hline
		\end{tabular}
		\caption{Output formats for the Fourier coefficients (s.v= singular values).}
		\label{table:formats}
	\end{table}

	\subsection{Fourier series}\label{ssec:fourier-series}
	
	Following our policy of efficiency, Fourier series are obtained using an expansion routine of the Fourier coefficients. We can suspect that Mathematica\texttrademark\, does not perform an expansion of the coefficients to obtain the Fourier series, but computes them one by one inside the \texttt{FourierSeries} function for two reasons. First, because Mathematica\texttrademark\, only works with truncated series and so a general expression for the coefficients is not needed. Second, because this approach avoids the problems associated to simplification of symbolic integrals using assumptions.\\
	
	 Nevertheless this strategy is inefficient because it does not take advantage of previously computed Fourier coefficients. As we have seen before, Mathematica\texttrademark\, can not obtain correctly the cosine coefficients of \( \cos^2 x \), but if we compute the first five terms of the cosine series of this function we get the right answer:
	\begin{verbatim}
		In[8]:=FourierCosSeries[Cos[x]^2,x,5]
		Out[8]:=1/2+1/2Cos[2x]
	\end{verbatim}
	
	In our case the upper limit of summation can be a positive integer or infinite. In the first case, a truncated series is returned; in the second, a symbolic series is displayed. In Table \ref{table:expansion} we summarize the syntax for the expansion routines along with the different Fourier series commands, see also Table \ref{table:formats} for comparison. \\
	\begin{table}
		\centering
		\begin{tabular}{l|l}
			Expansion routine & Series \\
			\hline\\\\[-1.9\medskipamount]
			\texttt{fouriercoeff\_expand(list of coeff,var,L,N)} & \texttt{fourier\_series(expr,var,L,N)} \\
			\texttt{cfouriercoeff\_expand(list of coeff,var,L,N)} & \texttt{cfourier\_series(expr,var,L,N)} \\
			\texttt{fouriersincoeff\_expand(list of coeff,var,L,N)} & \texttt{fouriersin\_series(expr,var,L,N)} \\
			\texttt{fouriercoscoeff\_expand(list of coeff,var,L,N)} & \texttt{fouriercos\_series(expr,var,L,N)} \\[7pt]
			\hline
		\end{tabular}
		\caption{Expansion routines and Fourier series syntax in \texttt{pdefourier}.}
		\label{table:expansion}
	\end{table}
	
	It is important to notice that displaying infinite series correctly has been a source of troubles in different CAS.
	This is mainly due to the fact that they are not evaluated, only displayed symbolically, so
	a simplification evident for a human might not be performed by a CAS. \\
	
	For illustration purposes, suppose that we want to compute the Fourier series of $h(x)=\sin 15 x$ on the interval $[-\pi,\pi]$.
	It is obvious that its Fourier series is exactly equal to $h(x)$, because its Fourier coefficients are 
	\[a_0=0, \, a_n=0, \, b_n=\delta_n^{15},\] 
	and we have:
	\[
	h(x) \sim \sum_{n=1}^\infty \delta_n^{15} \sin n x = \sin 15 x.
	\]
	
	However, infinite sums of expressions containing Kronecker delta functions can be hard to simplify to a
	single term, because it must be verified that the only index that does not vanish is indeed contained in the
	set of indices over which we are considering the sum. Integration routines in different CAS sometimes return
	the result for the Fourier integral formulas \eqref{eqn:fourier-coeff} in terms of some sort of Kronecker delta (Mathematica \textsuperscript{TM} returns it in terms of \texttt{DiscreteDelta}), which complicates the task of obtaining the Fourier series by an expansion
	of the coefficients. This is specially true if you are trying to implement Fourier series in a CAS or programming language that by itself doesn't have simplification algorithms of infinite series containing Kronecker deltas. \\
	
	Indeed, in the case of Maxima, even if the flag \texttt{simpsum} is set to true (by default it is set to false), some sums can be simplified but not those containing Kronecker deltas.
	\begin{verbatim}
	(%i3)	sum(1/n^2,n,1,inf);
	(%o3)	
	\end{verbatim}
	\[\sum_{n=1}^{\infty }{\left. \frac{1}{{{n}^{2}}}\right.}\]
	\begin{verbatim}
	(%i4)	simpsum:true$
	(%i5)	sum(1/n^2,n,1,inf);
	(%o5)	%pi^2/6
	(%i6)	block([n],declare(n,integer),
	assume(n>0),
	sum(kron_delta(15,n),n,1,inf));
	(%o6)	sum(kron_delta(15,n),n,1,inf)
	\end{verbatim}
	Following the approach described by algorithm \ref{algorithm:fcoeffs}, since we have already detected the set of indices of singular values of the coefficients, let's say \( \{k_1,k_2,\ldots, k_r\}\), whenever the user asks for the infinite Fourier series we can basically return
	\[
	\frac{1}{2}a_0+a_1\cos \frac{\pi x}{L}+b_1 \sin \frac{\pi x}{L}+\ldots+a_m\cos \frac{m\pi x}{L}+b_m \sin \frac{m\pi x}{L}+\sum_{n=m+1}^{\infty}\left( a_n \cos \frac{n \pi x}{L} +b_n \sin \frac{n \pi x}{L} \right),
	\]
	where \(m:=\max\{k_1,k_2,\ldots, k_r \}.\) In the next subsection we show an example of this situation in a Maple application. See also the examples of Fourier series in section \ref{sec-examples}.
	\subsection{Other approaches to the problem of symbolic computation of Fourier series}
	Given the importance of Fourier series, other people have developed different packages for their symbolic computation. We would like to give mention to three different approaches that we recently became aware of while writing this manuscript. All three approaches are developed as applications for Maple, an outstanding CAS with powerful algorithms. As we mentioned earlier in the introduction, we do not seek to compete with user packages or built-in routines of software such as Maple. Rather, we aim to show that even when working with such a mature and developed CAS, users have still needed to put in some effort to overcome all the previously exposed challenges symbolic Fourier series possess or are still present in a way. \\
	
	The first mention goes to the application \emph{FourierSeries} developed by Wilhelm Werner \cite{Werner}. This package was originally published in the year 2000 in the Maple Application Center, and later updated. The package handles what we call singular values exceptionally well. However, the package is mostly meant to compute Fourier coefficients, rather than working with Fourier series directly. Nonetheless, it offers some handy routines to visualize Fourier series. \\
	
\begin{figure}
		\centering
		\includegraphics[width=5in,height=2.4in]{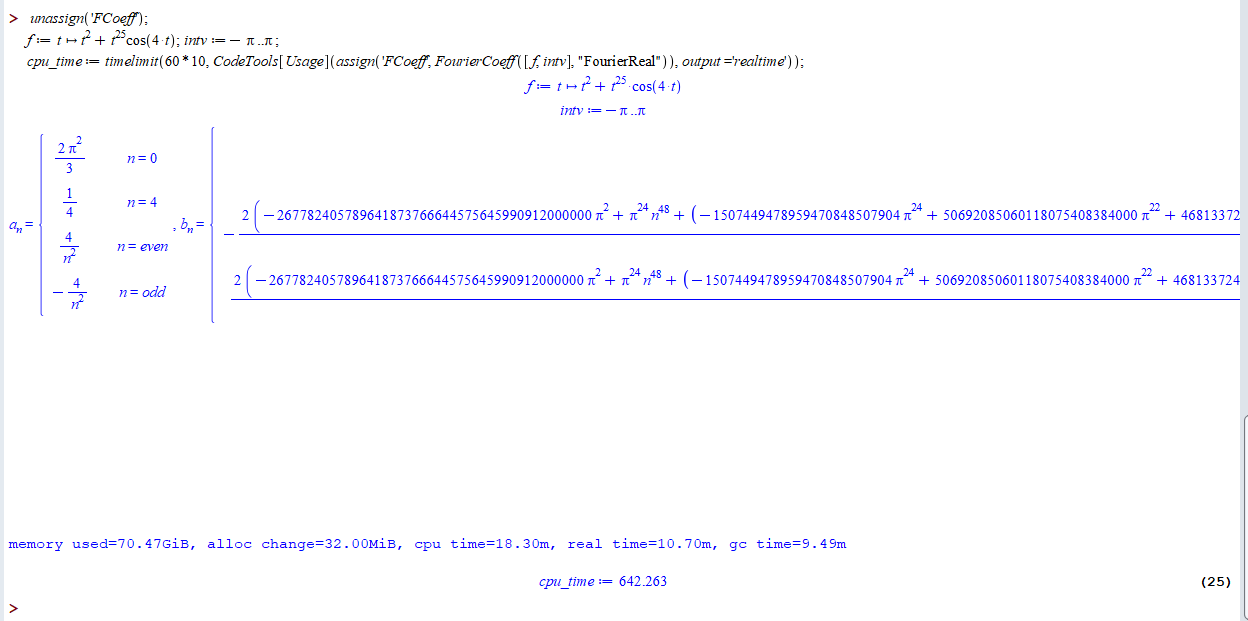}
		\caption{A computation of Fourier coefficients using Werner's Maple's user application \emph{FourierSeries}. In the case of a high-degree polynomial times a trigonometric function, the computation can take a huge amount of time.}
		\label{fig:werner}
	\end{figure}

	We would like to quote Werner's own words (appearing in the Application Demonstration worksheet) when talking about singular values in the coefficients:
	
	\begin{quotation}
		Probably many users encountered these problems and circumvented them by replacing symbolic integration with numerical integration methods. For graphical purposes this approach is acceptable, it is unsatisfactory in the context of computer algebra though. 
		Basically the approach which is taken by FourierSeries mimics hand calculation; the routine computes the generic solution and thereafter looks for critical cases. These exceptional cases are then handled separately. 
		The major problem of this approach is that the number of exceptional cases as well as their location is not known in advance. FourierSeries uses Maple's solve command to find exceptional cases and Maple's int command for integration;  the limitations of the package therefore are basically determined by the limitations of these Maple procedures, which, however, are known to be extremely powerful.
	\end{quotation}
	This is a clever approach, as it takes advantage of the Maple procedures which are well-known to be reliable. The strategy of using Maple's \texttt{solve} to find exceptional cases can have a big impact on computational time whenever we use a high-degree polynomial times a trigonometric function. Assuming the user has already installed Werner's package, setting a time limit of 600 seconds, computing the Fourier coefficients of \(f(t)=t^2+t^{25}\cos 4t\) returns a result although it crosses the time limit, see figure \ref{fig:werner}. Anyhow, this is a huge computational time. In section \ref{sec-examples}, we include the same example and its computational time when using \texttt{pdefourier}. \\
		
\begin{figure}
		\centering
		\includegraphics[width=5in,height=2.4in]{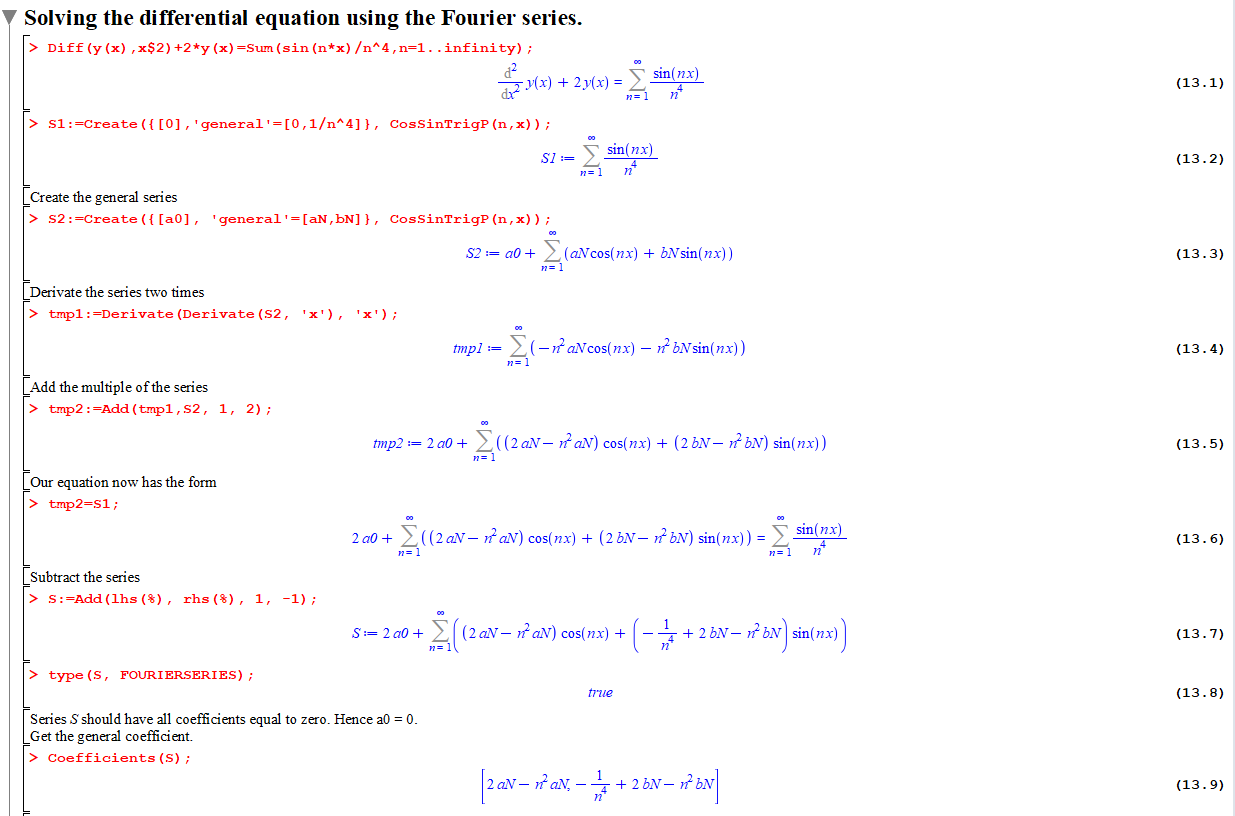}
		\caption{A screenshot showing the computational capabilities of Šrot's Maple package  \emph{FourierTrigSeries}, taken from the examples worksheet included in the package \cite{Srot}. A new data structure representation of the series allows the user to manipulate Fourier series. }
		\label{fig:karel}
	\end{figure}

	Our second mention goes to the application \emph{FourierTrigSeries}\footnote{This application was previously named \emph{FourierSeries}. Later, it was renamed to avoid confusion with Werner's application.} developed by Karel Šrot. The most updated version is available at Šrot's website \cite{Srot}. This package introduces a new data structure for the representation of trigonometric series which allows the user to manipulate Fourier series for many different computations. We recommend the reader to take a look at the example showing how to solve a differential equation involving Fourier series in the examples worksheet included in the package, see figure \ref{fig:karel}. This goes well-beyond the current capabilities of our implementation. Even though Šrot's handles trigonometric series with singular values pretty-well, trying the same example of a high-degree degree polynomial times a trigonometric function can take time as with Werner's package. \\
%
	
\begin{figure}
		\centering
		\includegraphics[width=5in,height=2.4in]{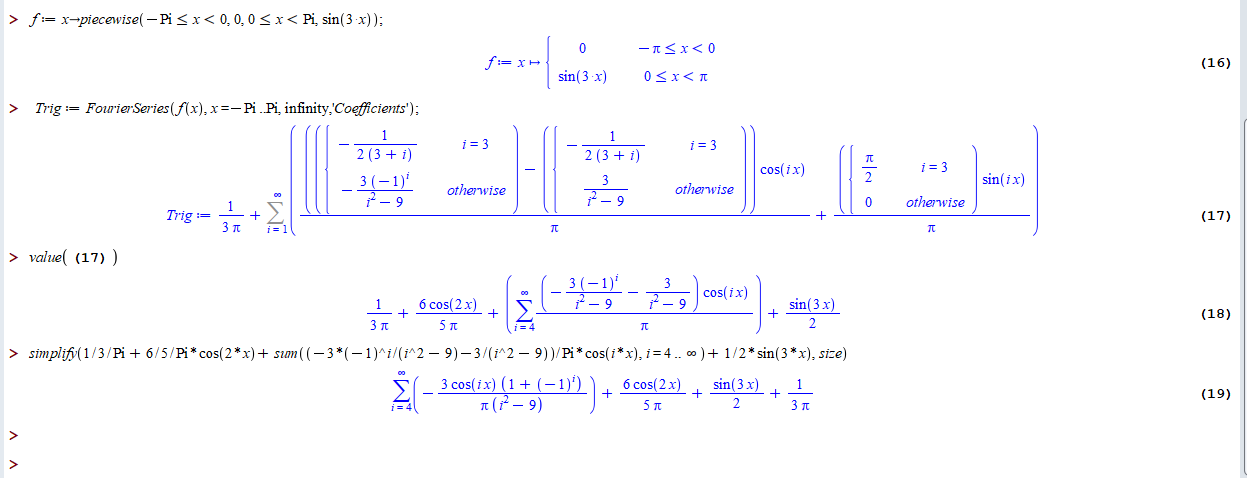}
		\caption{Screenshot showing how a direct expansion of coefficients with singular values is not necessarily simplified automatically. This example was generated using Maple's user application \emph{OrthogonalExpansions}. }
		\label{fig:orthogonal}
	\end{figure}

	Finally, there is the outstanding Maple application \emph{OrthogonalExpansions}, developed by Sergey Moiseev \cite{Moiseev}.  This package offers a collection of commands to compute one-dimensional and multi-dimensional orthogonal series expansions of a function. The expansions can be evaluated symbolically or numerically. This app includes routines for several orthogonal families, not only the trigonometric system. We would like to point out that the issues regarding expansion of Fourier series containing special cases for the coefficients (such as Kronecker deltas) are present in this package; however, thanks to the built-in capabilities of Maple, these issues are overcome easily, see figure \ref{fig:orthogonal}. Nevertheless, in our opinion this package should be considered as the standard in the topic of orthogonal expansions. Its approach paves the way to solve many other PDEs not considered here.
	
	\pagebreak
	\subsection{Numerical aspects of the 2D wave equation}
	The determination of Fourier coefficients in the case of the $2D$ wave equation requires computing the zeros of Bessel's
	functions. Unfortunately, Maxima 's numerical capabilities are rather limited. In particular, it does not have any built-in
	function for this task (although it has implemented the Bessel functions of first and second kind). Thus, it has been
	necessary to write a dedicated function for computing these zeros, and in the process of doing so several surprises
	regarding the implementation of Bessel functions in commercial computer algebra systems arose. \\
	
	Our solution for determining the zeros of $J_\nu(x)$ and $J'_\nu(x)$ is based on the papers \cite{GZ73} and \cite{IK91},
	which allows for a simple and efficient implementation of their algorithms, the main step being the computation of
	eigenvalues of certain matrices (we do this by importing the function \texttt{dgeev} from Netlib's \texttt{lapack}
	package). As we are interested exclusively in real zeros, we restrict the order in $J_\nu(x)$ to $\nu>-1$ (Lommel's theorem  \cite{Kerimov}),
	but otherwise we impose no restrictions on $\nu$. Thus, for example, the third zero of $J_{1/2}(x)$ is returned by
	the function \texttt{BesselJZeros} (mimicking the names in Maple\texttrademark\, and Mathematica\texttrademark\,):
	\begin{verbatim}
		(%i7) load(pdefourier)$
		(%i8)	BesselJZeros(1/2,3);
		(%o8)	9.42477796076938
	\end{verbatim}
	and the third zero of $J_{146225}(x)$ is
	\begin{verbatim}
		(%i9)	BesselJZeros(146225,3);
		(%o9)	146456.0070601201
	\end{verbatim}
	
	While Maple\texttrademark\, gives the result (even with better precision) effortlessly, to our surprise
	Mathematica\texttrademark\ could not calculate it (at least not in less than 10 minutes so we aborted the calculation), neither did the web-based engine Wolfram Alpha.
	Probably, this is due to the fact that they are using some variant
	of MacMahon formula, which fails for higher orders (see the comments in \cite{Wat18}).
	Other alternatives, such as Halley's method, also put a bound on the order, as it occurs in Matlab \cite{Nic20}. \\
	
	As mentioned, the algorithms in \cite{GZ73,IK91} can be used to generate the zeros of $J'_\nu(x)$. This is done
	by the function
	\texttt{BesselJdiffZeros}. The following table gives the first 5 zeros of the derivatives of $J'_\nu(x)$ for
	$0\leq\nu\leq 4$, and is to be compared with the one in \cite{Wei20}:
	\begin{verbatim}
		(%i10)	transpose(apply(matrix,makelist(BesselJdiffZeros(j,5,all),j,0,4)));
		(%o10)	matrix(
		[3.831705970207512,	1.84118378134066,	3.05423692822714,	4.201188941210528,	5.317553126083997],
		[7.015586669815619,	5.331442773525031,	6.706133194158456,	8.015236598375953,	9.28239628524161],
		[10.17346813506272,	8.536316366346284,	9.96946782308759,	11.345924310743,	12.68190844263889],
		[13.32369193631421,	11.70600490259207,	13.17037085601612,	14.58584828616704,	15.96410703773154],
		[16.47063005087759,	14.86358863390901,	16.34752231832178,	17.78874786606648,	19.19602880004888]
		)
	\end{verbatim}
	
	Up to our knowledge, Maple\texttrademark\, does not have a similar function for computing the zeros of derivatives of
	Bessel functions.

	\subsection{Solution of PDEs}
	In developing our library, we focused on solving the three main second-order linear partial differential equations with constant coefficients, namely, the heat, wave and Laplace equations.
	The mathematical details can be found in any standard textbook on PDEs or Fourier Analysis (for instance,
	see \cite{Hab13,Wal98}). Here, we just briefly discuss some technical aspects of the implementation, a list of examples will be given in the next section. \\
	
	It has been mentioned that Fourier series are obtained by performing an expansion of the Fourier coefficients. This approach is multipurpose. It is also useful to code the solution of the PDEs in an easier way. Although the package can solve the heat equation with a heat source $Q(x,t)$, we want to illustrate how the expansion sub-routines facilitate the process of finding solutions with a simpler instance of the general equation. Consider the following IBVP:
	\[
	\begin{cases}
		u_t = \kappa u_{xx} \qquad (x,t) \in [0,L] \times \mathbb{R^+} \\
		u(x,0)=f(x) \\
		u(0,t)=0 \\
		u(L,t)=0
	\end{cases}
	\]
	
	the solution is then given by
	\begin{align*}
		u(x,t) &= \sum_{n=1}^ \infty B_n \exp [-n^2 ( \frac{\pi \kappa}{L})^2 t ] \sin \frac{n \pi x}{L} \\
		B_n&=\frac{2}{L} \int_0^L f(x) \sin \frac{n \pi x }{L}
	\end{align*}
	
	Notice that for our purposes, it will be sufficient if we create a list of the form $[[C_n(t)],[ \, ] / [j,C_j(t)] ]$ where $C_n(t)=B_n \exp [-n^2 ( \frac{\pi \kappa}{L})^2 t ]$ and then use the expansion routine corresponding to a Fourier sine series on the interval $[0,L]$ and the space variable $x$. For instance, if we consider $f(x)=x^2(1-x),  \, L=1, \, \kappa=1$, a solution to the IBVP can be obtained in a few lines of code:
	
	\begin{verbatim}
		(%i11)	f(x):=if (0<=x and x<=1) then x^2*(1-x)$
		(%i12)	[[Bn],lsv]:fouriersincoeff(f(x),x,1)$
		(%i13)	Cn:Bn*exp(-n^2*%pi^2*t)$
		(%i14)	fouriersincoeff_expand([[Cn],[]],x,1,inf);
		(%o14)
	\end{verbatim}
	\[
	-\frac{4 \sum_{n=1}^{\infty }{\left. \frac{\left( 2 {{\left( -1\right) }^{n}}+1\right) \, {{\% e}^{-{{\ensuremath{\pi} }^{2}}\, {{n}^{2}} t}} \sin{\left( \ensuremath{\pi}  n x\right) }}{{{n}^{3}}}\right.}}{{{\ensuremath{\pi} }^{3}}}
	\]
	
	Here the list of singular values was empty, and no further work was required. The library
	command \texttt{mixed\_heat} automates all the process in the most general setting. \\

	Similar methods to obtain solutions for the three equations with different types of boundary conditions and domains have been implemented, dealing with general expressions (piecewise-defined or not) and taking care of the possible singular values in the coefficients.

	\section{Some examples}\label{sec-examples}
	
	\subsection{Fourier coefficients and series}
	
	Let $f(x)$ be as in \eqref{fexample}. Then we have the following:
	
	\begin{verbatim}
		(%i15)	fouriercoeff(3*x^2*cos(7*x),x,%pi);
		(%o15)
	\end{verbatim}
	\[
	[[-\frac{6}{49},\frac{12 \left( {{n}^{2}}+49\right) \, {{\left( -1\right) }^{n+1}}}{{{n}^{4}}-98 {{n}^{2}}+2401},0],[[7,\frac{98 {{\ensuremath{\pi} }^{2}}+3}{98},0]]]
	\]
	
	Indeed, see eq. \eqref{eqn:x2cos7x} and notice that $(n-7)^2(n+7)^2=n^4-98n^2+2401$. \\
	
	The case where Mathematica\texttrademark\,  fails to compute the cosine coefficients of equivalent trigonometric expressions \eqref{trigexample} is correctly solved by our package:
	
	\begin{verbatim}
		(%i16)	fouriercoscoeff(cos(x)^2,x,%pi);
		(%o16)
	\end{verbatim}
	\[
	[[\frac{1}{2},0],[[2,\frac{1}{2}]]]
	\]
	
	\begin{verbatim}
		(%i17)	fouriercoscoeff((1+cos(2*x))/2,x,%pi);
		(%o17)
	\end{verbatim}
	\[
	[[\frac{1}{2},0],[[2,\frac{1}{2}]]]
	\]
	Now, we consider an example of how Fourier series are displayed symbolically. If $g(x)=x^4$ on $[-\pi,\pi]$ we get the answer in the form used in textbooks:
	
	\begin{verbatim}
		(%i18)	fourier_series(x^4,x,%pi,inf);
		(%o18)
	\end{verbatim}
	\[
	8 \left( \sum_{n=1}^{\infty }{\left. \frac{\left( {{\ensuremath{\pi} }^{2}}\, {{n}^{2}}-6\right) \, {{\left( -1\right) }^{n}} \cos{\left( n x\right) }}{{{n}^{4}}}\right.}\right) +\frac{{{\ensuremath{\pi} }^{4}}}{5}
	\]
	Let us show an example of how to use the expansion routines to obtain the Fourier series, truncated or not, and how they handle singular values in the Fourier coefficients when displaying an infinite series. We show this example for the case of the sine system.
	
	\begin{verbatim}
		(%i19)	fcoeff:fouriersincoeff(x*cos(3*x),x,%pi);
		(fcoeff)
	\end{verbatim}
	\[
	[[\frac{2 n\, {{\left( -1\right) }^{n}}}{{{n}^{2}}-9}],[[3,-\frac{1}{6}]]]
	\]
	\begin{verbatim}
	(%i20)	fouriersincoeff_expand(fcoeff,x,%pi,5);
	(%o20)
	\end{verbatim}
	\[
	-\frac{5 \sin{\left( 5 x\right) }}{8}+\frac{8 \sin{\left( 4 x\right) }}{7}-\frac{\sin{\left( 3 x\right) }}{6}-\frac{4 \sin{\left( 2 x\right) }}{5}+\frac{\sin{(x)}}{4}
	\]
	\begin{verbatim}
	(%i21)	fouriersincoeff_expand(fcoeff,x,%pi,inf);
	(%o21)	
	\end{verbatim}
	\[2 \left( \sum_{n\mathop{=}4}^{\infty }{\left. \frac{n {{\left( \mathop{-}1\right) }^{n}} \sin{\left( n x\right) }}{{{n}^{2}}\mathop{-}9}\right.}\right) \mathop{-}\frac{\sin{\left( 3 x\right) }}{6}\mathop{-}\frac{4 \sin{\left( 2 x\right) }}{5}\mathop{+}\frac{\sin{(x)}}{4}\]
	
\begin{figure}
		\centering
		\includegraphics[trim=10pt 0 10pt 0cm, clip,scale=0.5]{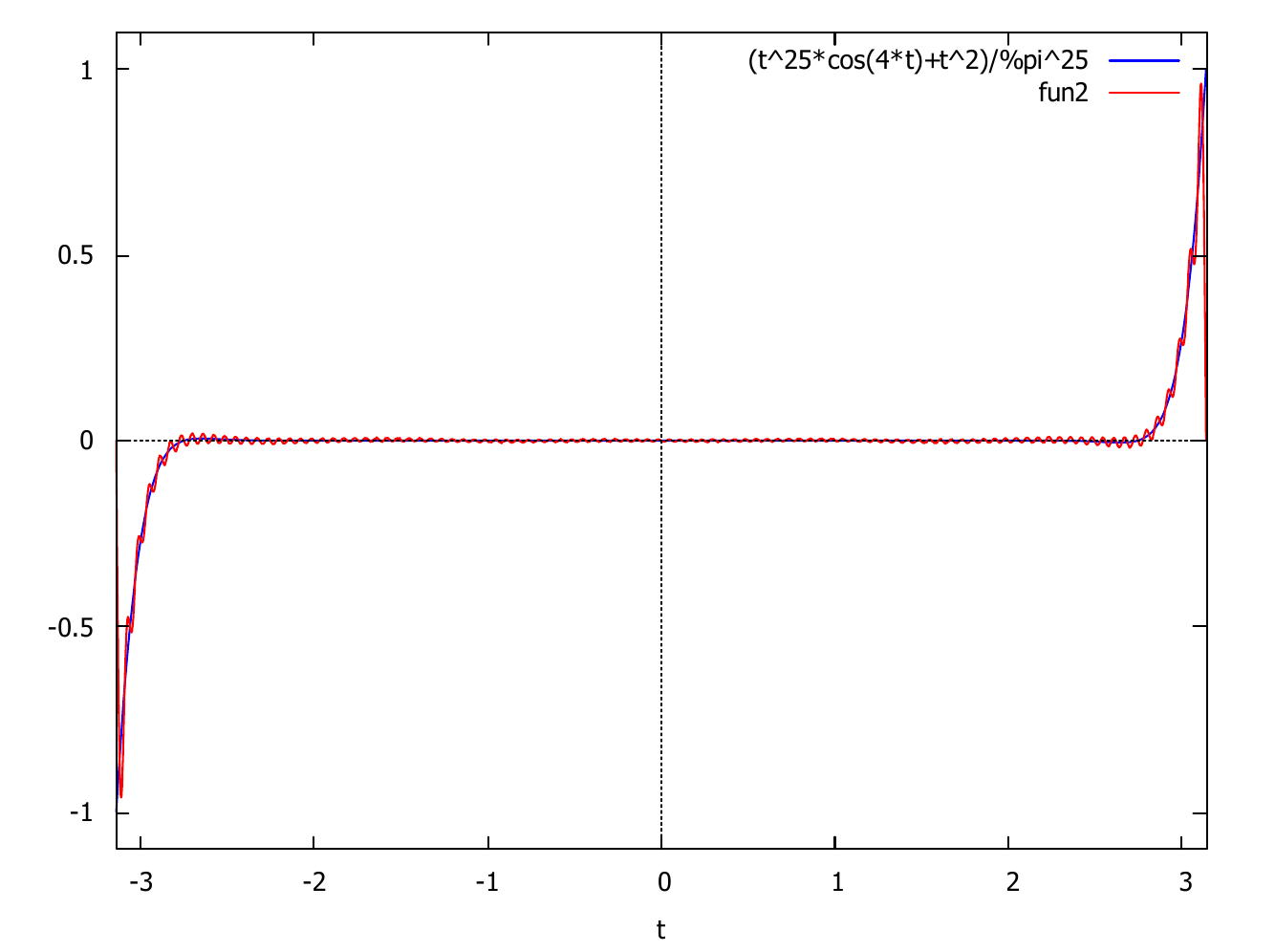}
		\caption{Plot of the Fourier series of \(t^2+t^{25}\cos 4t\). The scale factor \(1/\pi^{25}\) is included to help with floating-point rounding errors and have a nice-looking figure.}
		\label{fig:hdeg}
	\end{figure}

	In particular, following our approach as described in subsection \ref{ssec:fourier-series}, we avoid the issue of the evaluation of Kronecker delta functions inside an infinite sum:
	
	\begin{verbatim}
	(%i22)	fourier_series(sin(15*x),x,%pi,inf);
	(%o22) sin(15*x)
	\end{verbatim}
	
	Lastly, we offer an example of a piecewise-defined function having singular values. We will use the command \texttt{cfouriercoeff} to get the complex Fourier coefficients, so that we can compare our answer to the one returned by Mathematica\texttrademark.
	
	\begin{verbatim}
	(%i23)	f(x):=if x>= -%pi and x<0 then 0 elseif x>=0 and x<%pi then sin(3*x)$
	(%i24)	cfouriercoeff(f(x),x,%pi);
	(%o24)
	\end{verbatim}
	\[[[\frac{1}{3 \ensuremath{\pi} },-\frac{3 \left( {{\left( -1\right) }^{n}}+1\right) }{2 \ensuremath{\pi}  \left( {{n}^{2}}-9\right) }],[[3,-\frac{\% i}{4}]]]\]
	As we have shown in section \ref{ssec:hpw}, Mathematica \textsuperscript{TM} is not able to detect the singular value of the coefficient when $n=3$ for this piecewise-defined function. 
	\begin{rem}
		Notice that complex Fourier coefficients of real-valued functions satisfy the relation \(c_{-n}=\overline{c_{n}} \), where the bar denotes complex conjugation. Thus, it suffices to print the list of singular values with positive indices. The expansion routines handle the rest.
	\end{rem}
	To compare with Sergey Moiseev's \emph{OrthogonalExpansions} (see figure \ref{fig:orthogonal}), we can run the following command:
	\begin{verbatim}
		(%i25)	fourier_series(f(x),x,%pi,inf);
		(%o25)	
	\end{verbatim}

	\[\mathop{-}\left( \frac{3 \sum_{n\mathop{=}4}^{\infty }{\left. \frac{\left( {{\left( \mathop{-}1\right) }^{n}}\mathop{+}1\right)  \cos{\left( n x\right) }}{{{n}^{2}}\mathop{-}9}\right.}}{\ensuremath{\pi} }\right) \mathop{+}\frac{\sin{\left( 3 x\right) }}{2}\mathop{+}\frac{6 \cos{\left( 2 x\right) }}{5 \ensuremath{\pi} }\mathop{+}\frac{1}{3 \ensuremath{\pi} }\]
	Finally, we compare the time for computing the Fourier coefficients of a high-degree polynomial times a trigonometric function to compare with Werner's approach, see figure \ref{fig:werner}. Using Werner's package the computation takes more than 600 seconds of real time. With the pattern-matching approach it takes less than a second.
	\begin{verbatim}
		(%i26)	d:25$
		(%i27)	hdegexpr:(t^2+t^d*cos(4*t))$
		(%i28)	[[a0,an,bn],ls]:fouriercoeff(hdegexpr,t,%pi)$
		(%i29)	time(%);
		(%o29)	[0.21875]
	\end{verbatim}
	We can take a look at the leading term in the expression for \(b_n\):
	\begin{verbatim}
		(%i30)	factor(first(sum2list(expand(bn))));
		(%o30)	
	\end{verbatim}
	\[\mathop{-}\left( \frac{2 {{\left( \mathop{-}1\right) }^{n}} {{\ensuremath{\pi} }^{24}} {{n}^{49}}}{{{\left( n\mathop{-}4\right) }^{25}} {{\left( n\mathop{+}4\right) }^{25}}}\right) \]
	Since the coefficients have already been computed, we can expand them into an infinite series or a finite sum. Here, we use a finite sum to plot the result to compare with the original. When doing the plot, we divide by \(\pi^{25}\) to help with floating-point rounding errors. The plot can be seen in figure \ref{fig:hdeg}.
	\begin{verbatim}
	(%i31)	fcoeff:[[a0,an,bn],ls]$
	(%i32)	fs:fouriercoeff_expand(fcoeff,t,%pi,100)$
	(%i33)	time(%);
	(%o33)	[0.03125]
	(%i34)	wxplot2d([hdegexpr/%pi^d,expand(fs/%pi^d)],[t,-%pi,%pi]);
	(%t34)	 (Graphics) 
	\end{verbatim}
	\subsection{Partial differential equations}\label{secex}
	
	We will show some examples involving each type of equation (parabolic, hyperbolic and elliptic).
	Although \texttt{pdefourier} is able to solve these equations with many more boundary conditions, the ones
	given here are meant to give a panorama about how to use the package to solve PDEs. An excellent website that can be used as reference to consult many examples of Maple and Mathematica results when solving PDEs can be found in Nasser M. Abbasi's website \cite{A1200}.
	
	\begin{rem}
	In comparing results with Maple\texttrademark\, 2026 and Mathematica\texttrademark\ 14.2.1, we
	imposed a time limit of $600$ seconds (real time) to get an answer from both. If no answer was returned in such time, we considered it a failure. 
	\end{rem}
	
	\subsubsection{The linear second-order parabolic equation}
	Consider the general problem:
	\[
	\begin{cases}
		u_t=\kappa u_{xx} + v u_x + cu +Q(x,t) \mbox{ on } (x,t) \in [0,L] \times \mathbb{R}^+, \\
		u(x,0)=F(x), \\
		\alpha_1 u(0,t)+\beta_1 u_x(0,t) =h_1(t), \\
		\alpha_2 u(L,t)+\beta_2 u_x(L,t) =h_2(t).
	\end{cases}
	\]
	
	The syntax for solving it is
	\[
	\mbox{\texttt{mixed\_parabolic(Q(x,t),F(x),a1,b1,a2,b2,h1(t),h2(t),x,t,L,k,v,c,ord)}},
	\]
	where \texttt{ord} can be \texttt{inf}, for the complete series solution, or a natural number, for
	the truncated series.
	As an example, we will solve a problem with homogeneous boundary conditions
	\begin{equation}\label{eqn:mixed-heat}
	\begin{cases}
		u_t=\kappa u_{xx} - 9 u_x \quad (x,t) \in (0,1) \times \mathbb{R}^+, \\
		u(x,0)=e^{45x/10}(5\sin(\pi x)+9\sin(2\pi x)+2\sin(3\pi x)), \\
		u(0,t) = 0, \\
		u(1,t) =0.
	\end{cases}
	\end{equation}
	To avoid conflict with previously defined functions in previous sections, let us remove all definitions first and reload \texttt{pdefourier}.
	\begin{verbatim}
	(%i26)	kill(all);
	(%i1)	load(pdefourier)$
	\end{verbatim}
	To solve the problem above with \texttt{pdefourier} we would do
	\begin{verbatim}
		(%i2)	mixed_parabolic(0,exp(45*x/10)*(5*sin(%pi*x)+9*sin(2*%pi*x)+2*sin(3*%pi*x)),
		1,0,1,0,0,0,x,t,1,1,-9,0,inf);
		(%o2)
	\end{verbatim}
	\[{e^{\frac{9 x}{2}-\frac{81 t}{4}}}\, \left( 2 {e^{-9 {{\ensuremath{\pi} }^{2}} t}} \sin{\left( 3 \ensuremath{\pi}  x\right) }+9 {e^{-4 {{\ensuremath{\pi} }^{2}} t}} \sin{\left( 2 \ensuremath{\pi}  x\right) }+5 {e^{-{{\ensuremath{\pi} }^{2}} t}} \sin{\left( \ensuremath{\pi}  x\right) }\right) \]
	Both Maple\texttrademark\ and Mathematica\texttrademark\ solved it giving the same output as above. \\
	
\begin{figure}
		\centering
		\includegraphics[scale=0.66,trim={0 0 5cm 0},clip]{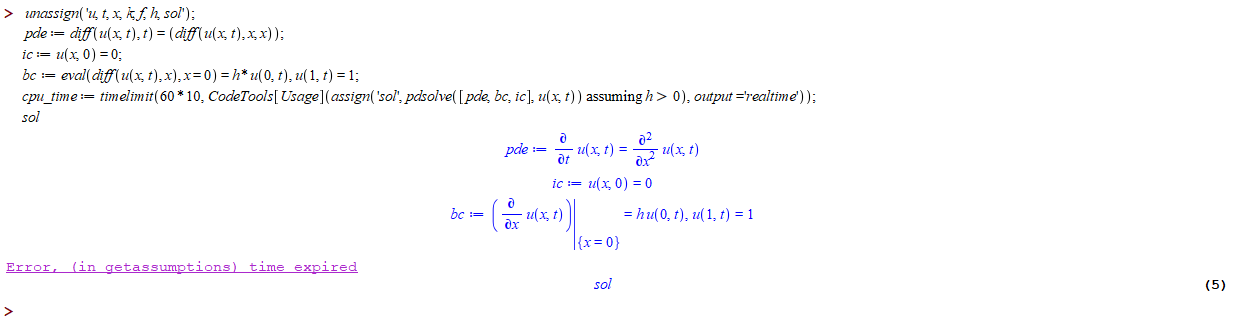}
		\caption{Screenshot showing that Maple can't solve problem \eqref{eqn:heat-loss} with the given time constrains.}
		\label{fig:maplepde}
	\end{figure}

	Now consider the following problem:
		\begin{equation}\label{eqn:heat-loss}
		\begin{cases}
			u_t=\kappa u_{xx}  \quad (x,t) \in (0,1) \times \mathbb{R}^+, \\
			u(x,0)=0, \\
			-hu(0,t)+u_x(0,t) = 0, \; h>0, \\
			u(1,t) =1.
		\end{cases}
	\end{equation}
	Both Maple\texttrademark\ and Mathematica\texttrademark\ fail to solve the following problem, see figure \ref{fig:maplepde}. 
	
	Our package solves it correctly
	\begin{verbatim}
	(%i3)	mixed_parabolic(0,0,-h,1,1,0,0,1,x,t,1,1,0,0,inf);
	(%03)   %lambda[n] are the solutions of 
	        -(h*sin(%lambda[n]))-%lambda[n]*cos(%lambda[n])=0
	\end{verbatim}
	\[\frac{h x\mathop{+}1}{h\mathop{+}1}\mathop{-}2 h\, \sum_{n\mathop{=}1}^{\infty }{\left. \frac{{{\% e}^{-\left( {{\lambda }_{n}^{2}} t\right) }} \sin{\left( {{\lambda }_n} \left( 1\mathop{-}x\right) \right) }}{{{\lambda }_n} \left( {{\cos{\left( {{\lambda }_n}\right) }}^{2}}\mathop{+}h\right) }\right.}\]
	\subsubsection{Laplace equation}
	Here we consider the Laplace equation on a wedge with Neumann conditions:
	\[
	\begin{cases}
		u_{rr}+\frac{1}{r} u_r + \frac{1}{r^2}u_{\theta \theta}=0 \mbox{ on } (r,\theta) \in [0,R] \times [0,\alpha] , \, 0<\alpha < 2\pi \\
		u(r,0)=0 \\
		u(r,\alpha)=0 \\
		u_r(R,\theta)=f(\theta)
	\end{cases}
	\]
	The syntax for Neumann problems in \texttt{pdefourier} is
	\[
	\mbox{\texttt{neumann\_laplace\_wedge(R,alpha,f(theta),theta,ord)}}.
	\]
	
	In this example we will consider $\alpha= \pi/2,R=1$, and an arbitrary
	function $f(\theta)$\footnote{The use of the percent sign is useful in some
		Maxima's graphical interfaces to get the symbol displayed as a Greek character.}:
	\begin{verbatim}
		(%i4)	neumann_laplace_wedge(1,%pi/2,f(%theta),%theta,inf);
		(%o4)
	\end{verbatim}
	\[
	\frac{2}{\pi} \sum_{n=1}^{\infty }{\left. \frac{\sin{\left( 2 \theta  n\right) } \int_{0}^{\frac{\ensuremath{\pi} }{2}}{\left. \operatorname{f}\left( \theta \right)  \sin{\left( 2 \theta  n\right) }d\theta \right.}\, {{r}^{2 n}}}{n}\right.}
	\]
	Mathematica is capable of solving it and it returns the same answer as above. However, if we use \(f(\theta)=\sin 2\theta\), it returns \( u(r,\theta)=0\) as an answer, which is not correct. 
	\begin{verbatim}
		In[14]:= 
		pde=D[u[r,theta],{r,2}]+D[u[r,theta],r]/r+D[u[r,theta],{theta,2}]/r^2==0;
		bc={Derivative[1,0][u][1,theta]==Sin[2*theta],u[r,Pi/2]==0,u[r,0]==0};
		sol=AbsoluteTiming[TimeConstrained[DSolve[{pde,bc},u[r,theta],{r,theta},
		Assumptions->{0<=r<=1&&0<=theta<=Pi/2}],60*10]];
		In[18]:= sol
		Out[18]= {1.35298,{{u[r,theta]->0}}}
	\end{verbatim}
	Our package is able to solve it correctly.
	\begin{verbatim}
		(%i5)	neumann_laplace_wedge(1,%pi/2,sin(2*%theta),%theta,inf);
		(%o5)	
	\end{verbatim}
	\[
	\frac{\sin(2\theta)r^2}{2},
	\]
	and so does Maple.
	\subsubsection{Wave equation}
	For this example, we will consider the wave equation with homogeneous boundary conditions and a driving term,
	a problem whose general expression is:
	\[
	\begin{cases}
		u_{tt}(x,t) = c^2 u_{xx}(x,t)+F(x,t) \quad (x,t) \in [0,L] \times \mathbb{R}^+\\
		u(x,0)=f(x) \\
		u_t(x,0)=g(x) \\
		u(0,t)=0 \quad u(L,t)=0
	\end{cases}
	\]
	
	The syntax for solving it using \texttt{pdefourier} is
	\[
	\mbox{\texttt{dirichlet\_wave(F(x,t),f(x),g(x),h\_1(t),h\_2(t),x,t,L,c,ord)}}.
	\]
	
	If $F(x,t)= \cos(\omega t) r(x)$ and $g(x)=0$, then:
	
	\begin{verbatim}
		(%i6)	assume(t>0,L>0)$
		(%i7)	dirichlet_wave(r(x)*cos(%omega*t),f(x),0,0,0,x,t,L,c,inf);
		(%o7)
	\end{verbatim}
	\begin{multline}
		\notag
		- \frac{\sum_{n=1}^\infty \sin ( \frac{\pi n x}{L }) \left( (2L^2 \cos ( \frac{\pi c n t }{L}) - 2L^2 \cos (\omega t) ) \int_0^L r(x) \sin (\frac{\pi n x}{L}) dx \right)}{\pi^2Lc^2 n^2 -L^3 \omega^2}  \\
		-\frac{ \sum_{n=1}^\infty ( 2L^2 \omega^2-2\pi^2c^2 n^2)\sin ( \frac{\pi n x}{L }) \cos ( \frac{\pi c n t }{L}) \int_0^L f(x) \sin (\frac{\pi n x}{L}) dx}{\pi^2Lc^2 n^2 -L^3 \omega^2}
	\end{multline}
	
	Since no assumptions were made about $\omega$, the solution corresponds to the case without resonance ($\omega \neq c n \pi/L$). This is Exercise 8.5.2(b) in \cite{Hab13}. For a different example, consider a vibrating clamped circular membrane,
	which is modeled in polar coordinates as
	\begin{equation}\label{eqn:wavedisk}
	\begin{cases}
		u_{tt}=c^2\left(u_{rr}+\frac{1}{r} u_r +\frac{1}{r^2}u_{\theta \theta}\right)\mbox{ on }(r,\theta,t)\in [0,R]\times [0,2\pi]
		\times [0,\infty [  \\
		u(r,\theta,0)=f(r,\theta) \\
		u_t(r,\theta,0)=g(r,\theta) \\
		u(R,\theta,t)=0 \\
		u(r,0,t)=u(r,2\pi,t) \\
		u_\theta (r,0,t)=u_\theta (r,2\pi,t) \
	\end{cases}
	\end{equation}
	
	The general case is handled by the function
	\[
	\mbox{\texttt{wave2d\_disk(c,R,f,g,k,l)},}
	\]
	where $k$ is the maximum order of the Bessel functions desired in the solution and $l$ is the number of positive zeros
	taken into account.
	In particular, for a circular membrane of unit radius, $c=1$, initial shape described by $f(r,\theta)=1-r^4$,
	and no initial velocities, we would issue the following commands:
	\begin{verbatim}
		(%i8)	f(r,theta):=1-r^4$
		(%i9)	g(r,theta):=0$
		(%i10)	wave2d_disk(1,1,f,g,3,2),numer$
	\end{verbatim}
	
	\begin{rem}
		There is a difference between this function and the remaining ones in the library. While we used expressions
		as arguments before (as \texttt{sin(5*\%pi*x)}), here we need to enter the names of declared functions (like \texttt{f,g}).
		Ultimately, this has to do with the fact that we are using numerical computations, but this is not an essential
		aspect, and as such will surely be modified in future versions to allow for expressions in the arguments.
	\end{rem}
	
	We have suppressed the output of the last command, since it contains a lot of terms with the aspect
	\[
	J_j(3,6{.}380161895923982r)\cos(6{.}380161895923982t) 2{.}7641436779037533\cdot 10^{-15}\sin(3\theta)\,.
	\]
	Obviously, these are numerical artifacts coming from rounding errors (as revealed by exponents like $10^{-15}$).
	In these cases, we can employ the function \texttt{chop} (provided in the package) analogous to
	the \texttt{Chop[]} encountered in Mathematica\texttrademark\,, to avoid rounding errors in the plot function.
	By default, \texttt{chop} discards terms of absolute value
	less than $10^{-12}$. An additional argument can be supplied, as in \texttt{chop(expr,exponent)}, to
	chop those terms whose absolute value is less that $10^{-exponent}$.
	\begin{verbatim}
		(%i11)	expr:chop(%,14);
		(expr)	1.3666632169857156*bessel_j(0.0,2.404825557695773*r)*cos(2.404825557695773*t)-
		0.4858370155994775*bessel_j(0.0,5.520078110286311*r)*cos(5.520078110286311*t)
	\end{verbatim}
	
\begin{figure}
		\centering
		\includegraphics[trim=10pt 0 10pt 1.5cm, clip,scale=0.5]{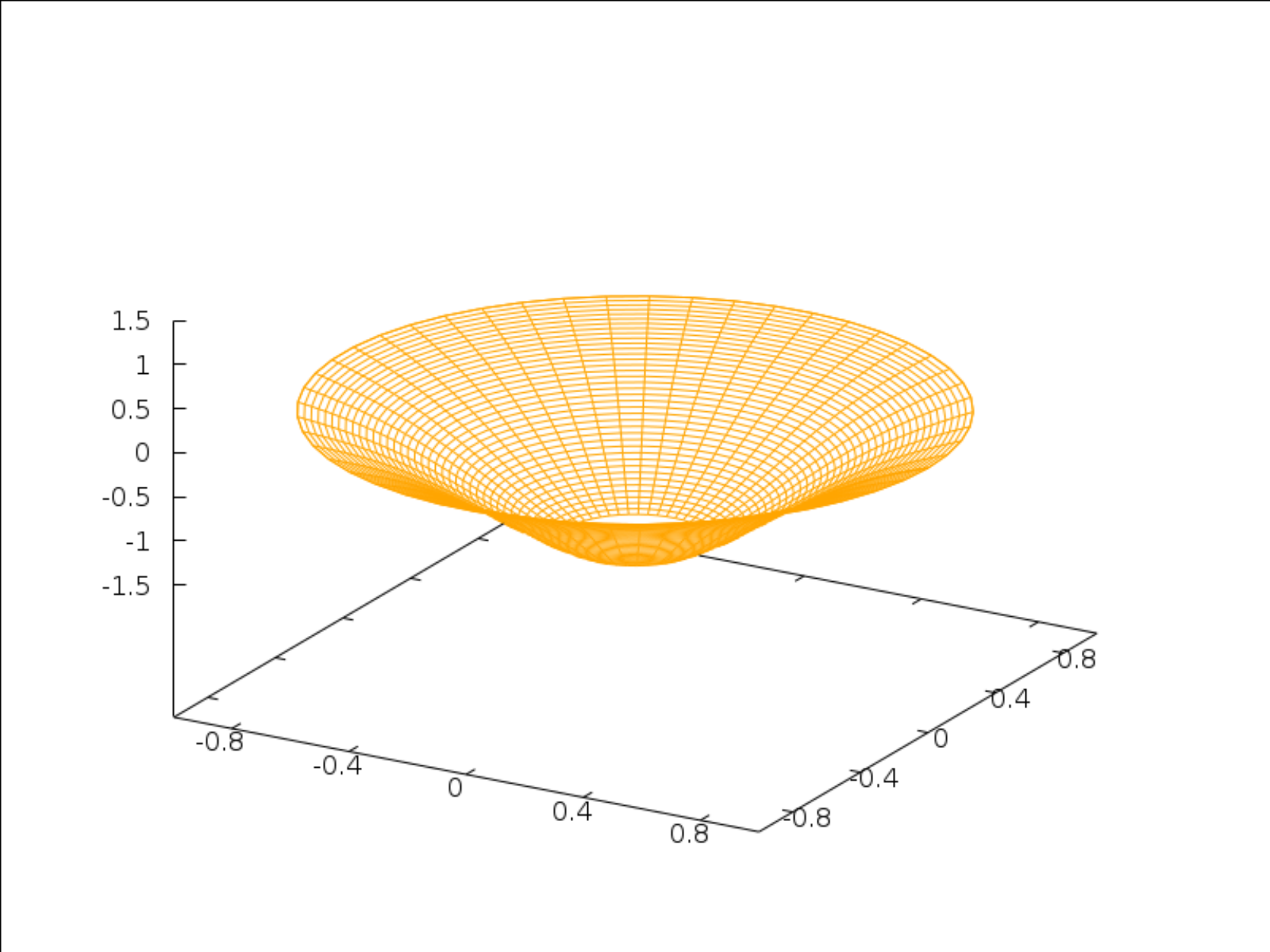}
		\caption{Capture of the animation generated by \texttt{wxanimate\_draw} for a  vibrating clamped circular membrane, modeled by the IBVP in \eqref{eqn:wavedisk} with \(f(r,\theta)=1-r^4,\) \(g(r,\theta)=0\), \(R=1\) and \(c=1\).}
		\label{fig:wave}
	\end{figure}

	The output has a format suitable for a graphical representation, or even to create an animation of the membrane
	motion\footnote{Here we use the command \texttt{wxanimate\_draw}, provided by the graphical frontend wxMaxima.}:
	\begin{verbatim}
		(%i12)	wxanimate_draw3d(s,makelist(i/10,i,0,24),
		surface_hide=true,zrange=[-2,2],
		color=orange,ytics=0.4,xtics=0.4,
		parametric_surface(r*cos(theta),r*sin(theta),subst(t=s,expr),r,0,1,theta,0,2*%pi)
		),wxanimate_framerate=6$
		(%t12)
	\end{verbatim}
	For a screencapture of the output see figure \ref{fig:wave}. When using Maxima versions \(\geq 5.48.1\), there is a bug in Gnuplot that prevents the option \texttt{surface\_hide=true} from working correctly. In that case, use \texttt{surface\_hide=false}.
	\section{Conclusions}
	We present a package for the Maxima CAS to compute Fourier series, and we use it to solve certain partial differential equations symbolically.  Throughout several examples, we show different difficulties that need to be addressed when trying to compute Fourier series symbolically. Even commercial computer algebra systems face different issues regarding the computation of Fourier series, or their application to the solution of certain partial differential equations. We describe a procedure using pattern-matching to catch special cases for the Fourier coefficients that we call singular values of the Fourier coefficients. We point out that our procedure allows us to compute Fourier series correctly and efficiently despite having certain limitations in integration routines, simplification of infinite series and other capabilities in our CAS. Our algorithm for computing Fourier series can be easily adapted to study partial differential equations in any other CAS, or in symbolic libraries such as SymPy. 

	\appendix
	\section{Installation of \texttt{pdefourier} } \label{app}
	
	The package \texttt{pdefourier} is available at \url{https://github.com/emmanuelroque/pdefourier}. Once the repository has
	been cloned or downloaded, the package can be installed by putting a copy of the files inside a folder contained in the
	environment variable \texttt{ file\_search\_maxima}. In a Linux box, such a system-wide location could be something like
	\texttt{/usr/share/maxima/5.42.2/share/contrib/}, while
	in a Windows environment typically it will be
	\begin{verbatim}
		C:\maxima-5.47.0/share/maxima/5.47.0/share/contrib
	\end{verbatim}
	(you may  need  administrator  rights  in  order  to  do  that in either case). The
	package can then be loaded  with the command \texttt{load(pdefourier)\$} once inside Maxima. Alternatively, users might look for the Maxima folder that gets created in their home folder (by default \texttt{C:\textbackslash Users\textbackslash myuser\textbackslash maxima}) and modify the file \texttt{maxima-init.mac} (sometimes you need to create it for the first time) with any text editor and add the following:
	\begin{verbatim}
		file_search_maxima: append(file_search_maxima,
		["C:/Users/myuser/my/path/to/pdefourier-folder/*.lisp",
		"C:/Users/myuser/my/path/to/pdefourier-folder/*.mac",
		"C:/Users/myuser/my/path/to/pdefourier-folder/*.wxm"])$
	\end{verbatim}

\end{document}